\documentclass[tikz,border=10pt]{article}

\usepackage[english]{babel}
\usepackage[utf8x]{inputenc}
\usepackage[T1]{fontenc}

\usepackage{tikz}
\usepackage{pgflibraryarrows}
\usepackage{pgflibrarysnakes}
\usetikzlibrary{decorations.markings}

\usepackage[a4paper,top=3cm,bottom=2cm,left=3cm,right=3cm,marginparwidth=1.75cm]{geometry}

\usepackage{amsmath}
\usepackage{graphicx}
\usepackage[colorinlistoftodos]{todonotes}
\usepackage[colorlinks=true, allcolors=blue]{hyperref} \usepackage{mathrsfs}
\usepackage{amssymb}
\usepackage{lineno}
\usepackage{bbding}
\usepackage{fancyhdr}\usepackage{graphics}
\usepackage{titlesec}\usepackage{titletoc}
\usepackage{indentfirst}
\usepackage{amsmath}\usepackage{amsfonts}\usepackage{latexsym}
\usepackage{eucal}\usepackage{epsfig}\usepackage{mathrsfs,amsmath,amsthm}
\usepackage{geometry}

\newtheorem{theorem}{Theorem}[section]
 \newtheorem{lemma}[theorem]{Lemma}
 \newtheorem{proposition}[theorem]{Proposition}

 \newcount\refno
\def\CC{{\mathbb C}}

 \def\RR{{\mathbb R}}
 \def\NN{{\mathbb N}}

 \def\SF{{\mathscr F}}
 \def\SH{{\mathscr H}}

 \title{\bf Paley type inequality on  Hardy spaces in the Dunkl setting
 \thanks{1} \footnote{E-mail:
huzhuoran010@163.com[ZhuoRan Hu].}}

\author{{ZhuoRan Hu}\\
{\small Department of Mathematics, Capital Normal University}\\
{\small Beijing 100048, China}}

\begin{document}

\maketitle \setcounter{page}{1} \pagestyle{myheadings}
 \markboth{Hu}{Paley type inequality on  Hardy space in the Dunkl setting}

\begin{abstract}
We investigate $\lambda$-Hilbert transform, $\lambda$-Possion integral and conjugate $\lambda$-Poisson integral on the atomic Hardy space in the Dunkl setting and establish a new version of  Paley type inequality which extends the results in \cite{F} and \cite{ZhongKai Li 3}.

\vskip .2in
 \noindent
 {\bf 2000 MS Classification:} 42B20, 42B25, 42A38.
 \vskip .2in
 \noindent
 {\bf 2000 Chinese Library Classification:} O178.
 \end{abstract}

\setcounter{page}{1}

\section{ Introduction and preliminaries }

For $0<p<\infty$, $L_{\lambda}^p(\RR)$ is the set of measurable functions satisfying
$ \|f\|_{L_{\lambda}^p}=\Big(c_{\lambda}\int_{\RR}|f(x)|^p|x|^{2\lambda}dx\Big)^{1/p}$ $<\infty$,
$c_{\lambda}^{-1}=2^{\lambda+1/2}\Gamma(\lambda+1/2)$,
and $p=\infty$ is the usual $L^\infty(\RR)$ space.
For $\lambda\geq0$, The Dunkl operator on the line is:
$$D_xf(x)=f'(x)+\frac{\lambda}{x}[f(x)-f(-x)]$$
involving a reflection part. The associated Fourier transfrom for the Dunkl setting for $f\in L_{\lambda}^1(\RR)$ is given by:
\begin{eqnarray}\label{fourier}
(\SF_{\lambda}f)(\xi)=c_{\lambda}\int_{\RR}f(x)E_\lambda(-ix\xi)|x|^{2\lambda}dx,\quad
\xi\in\RR , \,\,f\in L_{\lambda}^1(\RR).
\end{eqnarray}
$E_{\lambda}(-ix\xi)$ is the Dunkl kernel
$$E_{\lambda}(iz)=j_{\lambda-1/2}(z)+\frac{iz}{2\lambda+1}j_{\lambda+1/2}(z),\ \  z\in\CC$$
where $j_{\alpha}(z)$ is the normalized Bessel function
$$j_{\alpha}(z)=2^{\alpha}\Gamma(\alpha+1)\frac{J_{\alpha}(z)}{z^{\alpha}}=\Gamma(\alpha+1)\sum_{n=0}^{\infty}\frac{(-1)^n(z/2)^{2n}}{n!\Gamma(n+\alpha+1)} .$$
Since $j_{\lambda-1/2}(z)=\cos z$, $j_{\lambda+1/2}(z)=z^{-1}\sin z$, it follows that $E_0(iz)=e^{iz}$, and $\SF_{0}$ agrees with the usual Fourier transform. We assume $\lambda>0$ in
what follows.
And the associated $\lambda$-translation in Dunkl setting is
\begin{eqnarray}\label{tau}
 \tau_y f(x)=c_\lambda
     \int_{{\RR}}(\SF_{\lambda}f)(\xi)E(ix\xi)E(iy\xi)|\xi|^{2\lambda}d\xi,
     \ \ x,y\in{\RR} .
\end{eqnarray}
The $\lambda$-convolution$(f\ast_{\lambda}g)(x)$ of two appropriate functions $f$ and $g$ on $\RR$ associated to the $\lambda$-translation $\tau_t$ is defined by
$$(f\ast_{\lambda}g)(x)=c_{\lambda}\int_{\RR}f(t)\tau_{x}g(-t)|t|^{2\lambda}dt.$$
The "Laplace Equation" associated with the Dunkl setting is given by:
$$(\triangle_{\lambda}u)(x, y)=\left(D_x^2+ \partial_y^2\right) u(x, y)=\left(\partial_x^2+ \partial_y^2\right)u+ \frac{\lambda}{x}\partial_xu-\frac{\lambda}{x^2}\left(u(x, y)-u(-x, y)\right).$$
A $C^2$ function $u(x, y)$ satisfying $\triangle_{\lambda}u=0$ is  $\lambda$-harmonic.
When u and v are $\lambda$-harmonic functions satisfying $\lambda$-Cauchy-Riemann equations:
\begin{eqnarray}\label{a c r0}
\left\{\begin{array}{ll}
                                    D_xu-\partial_y v=0,&  \\
                                    \partial_y u +D_xv=0&
                                 \end{array}\right.
\end{eqnarray}
the function F(z)=F(x,y)=u(x,y)+iv(x,y)\,(z=x+iy)\, is a $\lambda$-analytic function.
 We define the Complex-Hardy spaces $H^p_\lambda(\RR^2_+)$ to be the set of
$\lambda$-analytic functions F=u+iv on $\RR^2_+$ satisfying
$$\|F\|_{H^p_\lambda(\RR^2_+)}=\sup\limits_{y>0}\left\{c_{\lambda}\int_{\RR}|F(x+iy)|^p|x|^{2\lambda}dx \right\}^{1/p}<\infty.$$

Throughout this paper, for any  $x_0\in\RR$, $\delta_0>0$, $I(x_0, \delta_0)$  is an Euclid interval with  center  $x_0$ and  radius $\delta_0$: $I(x_0, \delta_0)=\left\{x:|x-x_0|<\delta_0\right\}$ with $\delta_0<|{x_0}/2|$, and $I$ denotes as the set $I=\left(I(x_0,4\delta_0)\bigcup I(-x_0,4\delta_0)\right)^c.$

In this paper we define the atomic Hardy space in the Dunkl setting. A class of fundamental functions that we will call atoms will be introduced as following.  For $\frac{2\lambda}{2\lambda+1}< p\leq 1$, a Lebesgue measure  function $a(x)$ is  a $p_{\lambda}$-atom, if it satisfies the following conditions\\
(i) $\|a(x)\|_{L_{\lambda}^{\infty}}\lesssim \frac{1}{|I(x_0, \delta_0)|_{\lambda}^{1/p}}$, \\
(ii) ${\rm supp}\,a(x)\subseteq  I(x_0,\delta_0), with\ \  \delta_0<|{x_0}/2|$,  \\
(iii) $\int_{\RR} t^k a(t)|t|^{2\lambda}dt=0\ \ \ (k=0,1,2,3\ldots \kappa ), \ \kappa\geq2\left[(2\lambda+1)\frac{1-p}{p}\right]$.    \\

Let $\frac{2\lambda}{2\lambda+1}<  p\leq1$. Our  Hardy type space $H_{\lambda}^{p}(\RR)$ is constituted by all those $f\in S'(\RR)$ that can be represented by
\begin{eqnarray}\label{s}
f=\sum_{k}\lambda_{k}a_{k}(x),
\end{eqnarray}
being $\lambda_{j}\in \CC$ and $a_{j}$ is a $p_\lambda$-atom, for all $j\in\NN$, where $\sum_{j=0}^{\infty}|\lambda_{j}|^p<\infty$ and the series in (\ref{s}) converges in $S'(\RR)$.
We define on $H_{\lambda}^{p}(\RR)$ the norm $\|\cdot\|_{H_{\lambda}^{p}(\RR)}$ by
$$\|f\|_{H_{\lambda}^{p}(\RR)}=\inf\left(\sum_{k}|\lambda_{k}|^{p}\right)^{1/p},$$
where the infimum is taken over all those sequences $\{\lambda_j\}_{j\in\NN}\subset\CC$ such that $f$ is given by  (\ref{s}) for certain $p_{\lambda}$-atom $a_j$, $j\in\NN$.

We will establish the following version of Paley type inequality. Let $\frac{2\lambda}{2\lambda+1}<  p\leq1$, $k\ge p$ then there exists $C>0$ such that for every $f\in H_{\lambda}^{p}(\RR)$,
\begin{eqnarray*}
\int_0^{\infty}|(\SF_{\lambda}f)(\xi)|^k|\xi|^{(2\lambda+1)(k-1-k/p)+2\lambda}d\xi\le
c\|f\|_{H_{\lambda}^{p}(\RR)}^p.
\end{eqnarray*}
When $k=p$, we could obtain
$$\int_0^{\infty}|(\SF_{\lambda}f)(\xi)|^p|\xi|^{(2\lambda+1)(p-2)+2\lambda}d\xi\le
c\|f\|_{H_{\lambda}^{p}(\RR)}^p.$$
This  property is studied in \cite{R1,R2} for the classical case, in \cite{F} for the Dunkl setting on high dimension,  in \cite{ZhongKai Li 3} for the Dunkl setting on the upper half plane.

Our result is different to \cite{F,ZhongKai Li 3} that the  center of atoms in \cite{F} are fixed at the origin, and the functions in \cite{ZhongKai Li 3} are $\lambda$-analytic functions.

Throughout the paper we use the classic notation.  $S'(\RR)$ and $S(\RR)$ are the space of tempered distributions on $\RR$ and the Schwartz space on $\RR$ respectively.  $C$ will denote a  positive constant not necessary the same in each occurrence. We use $A\lesssim B$ to denote the estimate $|A|\leq CB$ for some  absolute universal constant $C>0$, which may vary from line to line,
$A\gtrsim B$ to denote the estimate $|A|\geq CB$ for some absolute universal constant $C>0$, $A\sim B$  to denote the estimate $|A|\leq C_1B$, $|A|\geq C_2B$ for some absolute universal constant $C_1, C_2$.

\section{Kernels and transforms associated with the Dunkl setting  }
The following properties are studied in \cite{ZhongKai Li 3}:
\begin{proposition}\cite{ZhongKai Li 3} \label{Poisson-a}
For $f\in L_{\lambda}^1(\RR)\bigcap L_{\lambda}^{\infty}(\RR)$, $x\in\RR, \ y\in(0,\infty)$, we can define $\lambda$-Hilbert transform, $\lambda$-Possion integral and conjugate $\lambda$-Poisson integral by
\begin{eqnarray*}
& &\SH_\lambda f(x)=c_{\lambda}\lim_{\epsilon\rightarrow0+}\int_{|t-x|>\epsilon} f(t)h(x,t)|t|^{2\lambda}dt, \\
& &(Pf)(x,y)=(f\ast_{\lambda}P_y)(x)=c_{\lambda}\int_{\RR}f(t)(\tau_xP_y)(-t)|t|^{2\lambda}dt, \\
& &(Qf)(x,y)=(f\ast_{\lambda}Q_y)(x)=c_{\lambda}\int_{\RR}f(t)(\tau_xQ_y)(-t)|t|^{2\lambda}dt,
\end{eqnarray*}
where
the $\lambda$-Hilbert kernel  is defined as :
$$h(x,t)=\frac{\lambda\Gamma(\lambda+1/2)}{2^{-\lambda-1/2}\pi}(x-t)\int_{-1}^1\frac{(1+s)(1-s^2)^{\lambda-1}
}{(x^2+t^2-2xts)^{\lambda+1}}ds,$$

$\lambda$-Poisson kernel $(\tau_xP_y)(-t)$ has the  representation
\begin{eqnarray}\label{D-Poisson-ker-11}
(\tau_xP_y)(-t)=
\frac{\lambda\Gamma(\lambda+1/2)}{2^{-\lambda-1/2}\pi}\int_0^\pi\frac{y(1+{\rm
sgn}(xt)\cos\theta)
}{\big(y^2+x^2+t^2-2|xt|\cos\theta\big)^{\lambda+1}}\sin^{2\lambda-1}\theta
d\theta,
\end{eqnarray}

and $(\tau_xQ_y)(-t)$ is the conjugate $\lambda$-Poisson kernel, with the following representation:
\begin{eqnarray}\label{D-conjugate-Poisson-ker-1}
(\tau_xQ_y)(-t)=
\frac{\lambda\Gamma(\lambda+1/2)}{2^{-\lambda-1/2}\pi}\int_0^\pi\frac{(x-t)(1+{\rm
sgn}(xt)\cos\theta)
}{\big(y^2+x^2+t^2-2|xt|\cos\theta\big)^{\lambda+1}}\sin^{2\lambda-1}\theta
d\theta.
\end{eqnarray}
$Pf(x, y)$ and  $Qf(x, y)$ satisfy the generalized Cauchy-Riemann system\,(\ref{a c r0})\,on $\RR_+^2$, and  $Pf(x, y)+iQf(x, y)$ is a $\lambda$-analytic function on the upper half plane.
\end{proposition}
The Hardy-Littlewood type inequalities for $F\in H_{\lambda}^p(\RR_+^2)$ are studied in \cite{ZhongKai Li 3}

\begin{theorem}\cite{ZhongKai Li 3}\label{us1}
Let $\frac{2\lambda}{2\lambda+1}<  p\leq1$, $k\ge p$ then there exists $C>0$ such that for every $f\in H_{\lambda}^{p}(\RR)$,
\begin{eqnarray*}
\int_0^{\infty}|(\SF_{\lambda}F)(\xi)|^k|\xi|^{(2\lambda+1)(k-1-k/p)+2\lambda}d\xi\le
c\|F\|_{H_{\lambda}^p(\RR_+^2)}^p.
\end{eqnarray*}
When $k=p$, we could obtain
$$\int_0^{\infty}|(\SF_{\lambda}F)(\xi)|^k|\xi|^{(2\lambda+1)(p-2)+2\lambda}d\xi\le
c\|F\|_{H_{\lambda}^p(\RR_+^2)}^p.$$
\end{theorem}

\begin{proposition}\cite{ZhongKai Li 3}\label{p5}If $\frac{2\lambda}{2\lambda+1}< p< l\le+\infty$,  $\delta=\frac{1}{p}-\frac{1}{l}$,
and $F(x, y)\in H^p_{\lambda}(\RR^2_+)$, $p\leq k<\infty$, then\\
{\rm(i)}
\begin{eqnarray}\label{i 1}
\left(\int_0^{+\infty} y^{k\delta(1+2\lambda)-1}\left(\int_{\RR}|F(x,y)|^{l}|x|^{2\lambda}dx \right)^{\frac{k}{l}}dy \right)^{\frac{1}{k}}\leq c\|F\|_{H_{\lambda}^p(\RR_+^2)}.
\end{eqnarray}
{\rm(ii)}
\begin{eqnarray}\label{i 2}
\left(\int_{\RR}|F(x,y)|^{l}|x|^{2\lambda}dx \right)^{\frac{1}{l}} \leq cy^{-(1/p-1/l)(1+2\lambda)}\|F\|_{H_{\lambda}^p(\RR_+^2)}.
\end{eqnarray}
When $p=\infty$, we could deduce the following from (\ref{i 2})
$$\sup_{x\in\RR}\big|F(x, y)\big| \leq cy^{-(1/p)(1+2\lambda)}\|F\|_{H_{\lambda}^p(\RR_+^2)} .$$
{\rm(iii)}
If $1\leq p<\infty$ and $F(x, y)=u(x, y)+iv(x, y)\in H^p_\lambda(\RR^2_+)$, then $F(x, y)$ is the $\lambda$-Poisson integrals of its boundary values $F(x)$, and $F(x)\in L^p_\lambda(\RR)$.
\end{proposition}

\begin{proposition} \label{estimate b}
$$\int_{I} \frac{1}{\left||x|-|x_0|\right|^{k}}dx\lesssim \delta_0^{1-k} \ \ \ (k>1).$$
\end{proposition}
\begin{proof}
\begin{eqnarray*}
\int_{I} \frac{1}{\left||x|-|x_0|\right|^{k}}dx\leq\int_{I(x_0,4\delta_0)^c} \frac{1}{\left||x|-|x_0|\right|^{k}}dx
=
\int_{I(0,4\delta_0)^c} \frac{1}{|x|^{k}}dx
\lesssim \delta_0^{1-k}.
\end{eqnarray*}
\end{proof}

\begin{proposition}\label{estimate a}
 $$\int_{-1}^{1}  \left(1-bs\right)^{-\lambda-1}(1+s)(1-s^2)^{\lambda-1}ds \leq C \frac{1}{1-|b|} ,\ \  \forall -1< b<1 \ ,\ \ \lambda>0$$
 C is depend on $\lambda$, and independent on $b$. ($C\thicksim 1/\lambda$)
\end{proposition}

\begin{proof}
CASE 1: when $0\leq b<1$.\\
It is obvious to see that when $0\leq b<1 \ ,\ \ \lambda>0$,
$$\int_{-1}^{0}  \left(1-bs\right)^{-\lambda-1}(1+s)(1-s^2)^{\lambda-1}ds \lesssim 1 .$$
By the formula of integration by parts and  $1-s\leq1-bs$ when $1\geq s\geq0$ ( $0\leq b<1 ,\ \ \lambda>0 $ ), we obtain:
\begin{eqnarray*}
\left|\int_{0}^{1} \left(1-bs\right)^{-\lambda-1}(1+s)(1-s^2)^{\lambda-1}ds\right|
 &\lesssim&
 \left|\int_{0}^{1}   \left(1-bs\right)^{-\lambda-1}(1-s)^{\lambda-1}ds\right|
\\&\lesssim&\frac{1}{\lambda}+\frac{\lambda+1}{\lambda}b\int_{0}^{1}\left(1-bs\right)^{-2}ds
\\&\lesssim&
\frac{1}{1-b}.
\end{eqnarray*}
Next we need to prove when $-1<b\leq0$
$$\int_{-1}^{1}  \left(1-bs\right)^{-\lambda-1}(1+s)(1-s^2)^{\lambda-1}ds \lesssim \frac{1}{1+b}.$$
CASE 2: when $-1<b\leq0$.\\
Obviously the following inequality holds:
$$\int_{0}^{1}  \left(1-bs\right)^{-\lambda-1}(1+s)(1-s^2)^{\lambda-1}ds \lesssim 1 .$$
By the formula of integration by parts and  $1+s\leq1-bs$ when $-1\leq s\leq0$ ( $-1<b\leq0 ,\ \ \lambda>0 $ ), we obtain:
\begin{eqnarray*}
 \left|\int_{-1}^{0} \left(1-bs\right)^{-\lambda-1}(1+s)(1-s^2)^{\lambda-1}ds\right|
 &\lesssim&
 \left|\int_{-1}^{0}   \left(1-bs\right)^{-\lambda-1}(1+s)^{\lambda}ds\right|
\\&\lesssim&\frac{1}{\lambda+1}-b\int_{-1}^{0}\left(1-bs\right)^{-1}ds
\\&\lesssim&
-\ln(1+b)
\\&\lesssim&
\frac{1}{1+b}.
\end{eqnarray*}
By CASE 1 and CASE 2, the inequality:
$$\int_{-1}^{1}  \left(1-bs\right)^{-\lambda-1}(1+s)(1-s^2)^{\lambda-1}ds \leq C \frac{1}{1-|b|} ,\ \  \forall -1< b<1 \ ,\ \ \lambda>0$$
holds. Hence the proposition holds.
\end{proof}
Thus we could obtain the following Proposition\,\ref{estimate c} and Proposition\,\ref{estimate a}
\begin{proposition}\label{estimate c}
 $$\int_{-1}^{1}  \left(1-bs\right)^{-\lambda-1}(1-s)(1-s^2)^{\lambda-1}ds \leq C \frac{1}{1-|b|} ,\ \  \forall -1< b<1 \ ,\ \ \lambda>0$$
 C is depend on $\lambda$, and independent on $b$.
\end{proposition}
The following Proposition\,\ref{estimate d} could be obtained in a way similar to Proposition\,\ref{estimate a}
\begin{proposition}\label{estimate d}
 $$\int_{-1}^{1}  \left(1-bs\right)^{-\lambda-1}(1-s^2)^{\lambda-1/2}ds \leq C \frac{1}{1-|b|} ,\ \  \forall -1< b<1 \ ,\ \ \lambda>0$$
 C is depend on $\lambda$, and independent on $b$. (In fact $C\thicksim \frac{2}{2\lambda+1}$)
\end{proposition}

\begin{theorem}[$\lambda$-Hilbert transform]\label{p atom}
For $\frac{2\lambda}{2\lambda+1}< p\leq 1$, if  $a(t)$ is  a $p_{\lambda}$-atom,  with vanishing order $\displaystyle{\kappa\geq2\left[(2\lambda+1)\frac{1-p}{p}\right]}$ then the following  holds:
\begin{eqnarray*}
 \int_{\RR}  |\SH_\lambda a(x)|^p |x|^{2\lambda}dx \leq C,
\end{eqnarray*}
C is depend on $\lambda$ and p.
\end{theorem}

\begin{proof}
Assume first that $x_0>0.$ Let $\kappa=2\left[(2\lambda+1)\frac{1-p}{p}\right].$
Thus $\kappa$ is an even integer. Let $n=\kappa/2$.
We could write the above integral as:
\begin{eqnarray*}
\int_{\RR}  |\SH_\lambda a(x)|^p |x|^{2\lambda}dx
&=&\int_{I(x_0,4\delta_0)\bigcup I(-x_0,4\delta_0)}  |\SH_\lambda a(x)|^p |x|^{2\lambda}dx    \\
&+& \int_{\left(I(x_0,4\delta_0)\bigcup I(-x_0,4\delta_0)\right)^c}  |\SH_\lambda a(x)|^p |x|^{2\lambda}dx \\
&=&I+II.
\end{eqnarray*}

We could see the following inequality holds:
\begin{eqnarray}\label{9}
4^{2\lambda+1}\int_{x_0-\delta_0}^{x_0+\delta_0}|x|^{2\lambda}dx
=\int_{4x_0-4\delta_0}^{4x_0+4\delta_0} |x|^{2\lambda}dx
\geq\int_{4x_0-\delta_0}^{4x_0+\delta_0}|x|^{2\lambda}dx.
\end{eqnarray}
By \cite{ZhongKai Li 3}\,[Theorem 5.7], together with (\ref{9}) we  obtain:
\begin{eqnarray}\label{8}
I&=&\int_{I(x_0,4\delta_0)\bigcup I(-x_0,4\delta_0)}  |\SH_\lambda a(x)|^p |x|^{2\lambda}dx  \nonumber \\
&\leq& \left(\int_{I(x_0,4\delta_0)\bigcup I(-x_0,4\delta_0)}  |\SH_\lambda a(x)|^2 |x|^{2\lambda}dx\right)^{p/2}
\left(\int_{I(x_0,4\delta_0)\bigcup I(-x_0,4\delta_0)}  |x|^{2\lambda}dx\right)^{1-p/2} \nonumber \\
&\leq& C.
\end{eqnarray}

Next we need to prove:
\begin{eqnarray}
II=\int_{(I(x_0,4\delta_0)\bigcup I(-x_0,4\delta_0))^c}  |\SH_\lambda a(x)|^p |x|^{2\lambda}dx \leq C.
\end{eqnarray}

By Proposition\,\ref{Poisson-a}, when $x\in \left(I(x_0,4\delta_0)\bigcup I(-x_0,4\delta_0)\right)^c$ we could write $\SH_\lambda a(x)$ as:
$$\SH_\lambda a(x)=c_{\lambda}\int a(t) h(x, t)|t|^{2\lambda}dt .$$
Next we need to estimate $\SH_\lambda a(x)$ when $x\in I=\left(I(x_0,4\delta_0)\bigcup I(-x_0,4\delta_0)\right)^c$    ($\delta_0<|x_0/2|$).

Notice that $t\in supp\,a(t)\subseteq I(x_0, \delta_0)$.  When $x\geq0,\  or \  x<-2x_0$, the following inequality
\begin{eqnarray}\label{8q}
|x-x_0|\lesssim\left(\langle x,x_0\rangle_s\right)^{1/2}
\end{eqnarray}
holds.
It is also obvious to see that the following inequalities hold:
\begin{eqnarray}\label{8q1}
|x-x_0s|\lesssim\left(\langle x,x_0\rangle_s\right)^{1/2},
\end{eqnarray}
\begin{eqnarray}\label{8q2}
|xs-x_0|\lesssim\left(\langle x,x_0\rangle_s\right)^{1/2},
\end{eqnarray}
\begin{eqnarray}\label{8qq}
|x_0+t-2xs|&\leq&|x_0-xs|+|t-xs|\\ \nonumber
&\leq&\left(\langle x,x_0\rangle_s\right)^{1/2}+\left(\langle x,x_0\rangle_s\right)^{1/2}+|t-x_0| \\ \nonumber
&\leq&3\left(\langle x,x_0\rangle_s\right)^{1/2}.\nonumber
\end{eqnarray}

From (\ref{8qq}), we could obtain the following inequality:
\begin{eqnarray}\label{taylor}
|\delta_1|\leq3|t-x_0|\left(\langle x,x_0\rangle_s\right)^{1/2}.
\end{eqnarray}

For $\delta_0<|x_0/2|, \ and\ x\in I=\left(I(x_0,4\delta_0)\bigcup I(-x_0,4\delta_0)\right)^c,$ from (\ref{taylor}), we could have:
\begin{eqnarray}\label{1}
\left| \frac{\delta_1}{\langle x,x_0\rangle_s} \right| \leq \left| \frac{3|t-x_0|}{\left(\langle x,x_0\rangle_s\right)^{1/2}} \right|\leq \left| \frac{3|t-x_0|}{\left||x|-|x_0|\right|} \right|\leq \left|\frac{3\delta_0}{4\delta_0} \right|=3/4.
\end{eqnarray}

We could see that:
\begin{eqnarray}\label{21}
\frac{x-t}{\langle x,t\rangle_s^{\lambda+1}}&=&\frac{x-x_0}{\langle x,t\rangle_s^{\lambda+1}}+\frac{x_0-t}{\langle x,t\rangle_s^{\lambda+1}}
\\&=&\nonumber A+B
\end{eqnarray}

By  Taylor expansion of  $\left(1+\frac{\delta_1}{\langle x,x_0\rangle_s}\right)^{-\lambda-1}$, when $x\in \left[-2x_0, 0\right]^c\bigcap \left(I(x_0,4\delta_0)\bigcup I(-x_0,4\delta_0)\right)^c $, we could obtain:
\begin{eqnarray}\label{22}
A&=&\frac{x-x_0}{\langle x,x_0\rangle_s^{\lambda+1} \left(1+\frac{\delta_1}{\langle x,x_0\rangle_s}\right)^{\lambda+1}}
\\&=&\frac{x-x_0}{\langle x,x_0\rangle_s^{\lambda+1}}\Bigg[1+\frac{\lambda+1}{1}\left(\frac{-\delta_1}{\langle x,x_0\rangle_s}\right)^1 \nonumber +\frac{(\lambda+1)(\lambda+2)}{2!}\left(\frac{-\delta_1}{\langle x,x_0\rangle_s}\right)^2+\cdots\\&+&\frac{(\lambda+1)_{n}}{(n)!}\left(\frac{-\delta_1}{\langle \nonumber x,x_0\rangle_s}\right)^{n}+\frac{(\lambda+1)_{n+1}}{(n+1)!}\left(\frac{1}{1+\xi_1}\right)^{\lambda+n+2}\left(\frac{-\delta_1}{\langle x,x_0\rangle_s}\right)^{n+1} \Bigg],\nonumber
\end{eqnarray}
and
\begin{eqnarray}\label{23}
B&=&\frac{x_0-t}{\langle x,x_0\rangle_s^{\lambda+1} \left(1+\frac{\delta_1}{\langle x,x_0\rangle_s}\right)^{\lambda+1}}
\\&+&\frac{x_0-t}{\langle x,x_0\rangle_s^{\lambda+1}}\Bigg[1+\frac{\lambda+1}{1}\left(\frac{-\delta_1}{\langle x,x_0\rangle_s}\right)^1 \nonumber +\frac{(\lambda+1)(\lambda+2)}{2!}\left(\frac{-\delta_1}{\langle x,x_0\rangle_s}\right)^2+\cdots\\&+&\frac{(\lambda+1)_{n-1}}{(n-1)!}\left(\frac{-\delta_1}{\langle \nonumber x,x_0\rangle_s}\right)^{n-1}+\frac{(\lambda+1)_{n}}{(n)!}\left(\frac{1}{1+\xi_2}\right)^{\lambda+n+1}\left(\frac{-\delta_1}{\langle x,x_0\rangle_s}\right)^{n} \Bigg].
\end{eqnarray}

By (\ref{1}),  $\xi_1,  \xi_2 \in[-3/4,3/4]$.
Thus the inequality holds:
\begin{eqnarray}\label{3}
\left(\frac{1}{1+\xi_1}\right)\leq \left(\frac{1}{1-3/4}\right)\leq 4 \ and\  \left(\frac{1}{1+\xi_2}\right)\leq \left(\frac{1}{1-3/4}\right)\leq 4.
\end{eqnarray}

Thus from Proposition\,\ref{Poisson-a}, Formulas\,(\ref{1}, \ref{21}, \ref{22}, \ref{23}, \ref{8q}, \ref{8qq}, \ref{taylor}), together with the vanishing property of $a(t)$, we  obtain:
for $x\in \left[-2x_0, 0\right]^c\bigcap \left(I(x_0,4\delta_0)\bigcup I(-x_0,4\delta_0)\right)^c $
\begin{eqnarray}\label{poisson inequality 1}
\left|\SH_\lambda a(x)\right|
&\lesssim&\left|\int\int_{-1}^1 |a(t)|\left|\frac{x-x_0}{\langle x,x_0\rangle_s^{\lambda+1}}\left(\frac{-\delta_1}{\langle x,x_0\rangle_s}\right)^{n+1}\right|(1+s)(1-s^2)^{\lambda-1}ds|t|^{2\lambda}dt\right|\nonumber\\&+&
\int\int_{-1}^1 |a(t)|\left|\frac{x_0-t}{\langle x,x_0\rangle_s^{\lambda+1}}\left(\frac{-\delta_1}{\langle x,x_0\rangle_s}\right)^{n}\right|(1+s)(1-s^2)^{\lambda-1}ds|t|^{2\lambda}dt\nonumber\\
&\lesssim&  |I(x_{0} ,\delta_{0})|_{\lambda}^{1-(1/p)} \int_{-1}^{1} \frac{(\delta_0)^{n+1}}{\left(\langle x,x_0\rangle_s\right)^{n/2+\lambda+1}} (1+s)(1-s^2)^{\lambda-1}ds.
\end{eqnarray}

It is easy to see the following inequality holds:
\begin{eqnarray}\label{4}
 \int_{-1}^{1} \frac{(\delta_0)^{n+1}(1+s)(1-s^2)^{\lambda-1}}{\left(\langle x,x_0\rangle_s\right)^{n/2+\lambda+1}} ds
 \leq
 \int_{-1}^{1} \frac{(\delta_0)^{n+1}(1+s)(1-s^2)^{\lambda-1}}{\left||x|-|x_0|\right|^n \left(\langle x,x_0\rangle_s\right)^{\lambda+1}} ds.
\end{eqnarray}

By Proposition(\ref{estimate a}), we obtain the following
\begin{eqnarray}\label{5}
\int_{-1}^1\frac{(1+s)(1-s^2)^{\lambda-1}
}{(x^2+x_0^2-2xx_0s)^{\lambda+1}}ds
&=&
\left( \frac{1}{x^2+x_0^2}\right)^{\lambda+1}\int_{-1}^1(1+s)(1-s^2)^{\lambda-1}
\left(1-\frac{2xx_0 s}{x^2+x_0^2 }\right)^{-\lambda-1}ds  \nonumber\\
&\leq&
C\frac{1}{(x^2+x_0^2)^{\lambda} (|x|-|x_0|)^2}.
\end{eqnarray}
Thus from Formulas\,(\ref{4}, \ref{5}) we have
\begin{eqnarray}\label{101}
& &|I(x_{0} ,\delta_{0})|_{\lambda}^{1-(1/p)} \int_{-1}^{1} \frac{(\delta_0)^{n+1}}{\left(\langle x,x_0\rangle_s\right)^{n/2+\lambda+1}} (1+s)(1-s^2)^{\lambda-1}ds
\nonumber\\&\lesssim&
|I(x_{0} ,\delta_{0})|_{\lambda}^{1-(1/p)} \frac{(\delta_0)^{n+1}}{\left||x|-|x_0|\right|^{n+2}\left||x|+|x_0|\right|^{2\lambda}}.
\end{eqnarray}
When $x\in \left[-2x_0, 0\right]^c\bigcap \left(I(x_0,4\delta_0)\bigcup I(-x_0,4\delta_0)\right)^c $, we obtain:
\begin{eqnarray}\label{104}
\left|\SH_\lambda a(x)\right|\leq C |I(x_{0} ,\delta_{0})|_{\lambda}^{1-(1/p)} \frac{(\delta_0)^{n+1}}{\left||x|-|x_0|\right|^{n+2}\left||x|+|x_0|\right|^{2\lambda}}.
\end{eqnarray}

Then  we need to consider the case  when $x\in \left[-2x_0, 0\right]\bigcap \left(I(x_0,4\delta_0)\bigcup I(-x_0,4\delta_0)\right)^c $. We could see that
\begin{eqnarray}\label{31}
\frac{x-t}{\langle x,t\rangle_s^{\lambda+1}}&=&\frac{x-x_0s}{\langle x,t\rangle_s^{\lambda+1}}+\frac{x_0-t}{\langle x,t\rangle_s^{\lambda+1}}+\frac{x_0(s-1)}{\langle x,x_0\rangle_s^{\lambda+1}}
\\&=&\nonumber C+D+E.
\end{eqnarray}
By the Taylor expansion, we could obtain:
\begin{eqnarray}\label{32}
C&=&\frac{x-x_0s}{\langle x,x_0\rangle_s^{\lambda+1} \left(1+\frac{\delta_1}{\langle x,x_0\rangle_s}\right)^{\lambda+1}}
\\&=&\frac{x-x_0s}{\langle x,x_0\rangle_s^{\lambda+1}}\Bigg[1+\frac{\lambda+1}{1}\left(\frac{-\delta_1}{\langle x,x_0\rangle_s}\right)^1 \nonumber +\frac{(\lambda+1)(\lambda+2)}{2!}\left(\frac{-\delta_1}{\langle x,x_0\rangle_s}\right)^2+\cdots\\&+&\frac{(\lambda+1)_{n}}{(n)!}\left(\frac{-\delta_1}{\langle \nonumber x,x_0\rangle_s}\right)^{n}+\frac{(\lambda+1)_{n+1}}{(n+1)!}\left(\frac{1}{1+\xi_3}\right)^{\lambda+n+2}\left(\frac{-\delta_1}{\langle x,x_0\rangle_s}\right)^{n+1} \Bigg],\nonumber
\end{eqnarray}

\begin{eqnarray}\label{33}
D&=&\frac{x_0-t}{\langle x,x_0\rangle_s^{\lambda+1} \left(1+\frac{\delta_1}{\langle x,x_0\rangle_s}\right)^{\lambda+1}}
\\&+&\frac{x_0-t}{\langle x,x_0\rangle_s^{\lambda+1}}\Bigg[1+\frac{\lambda+1}{1}\left(\frac{-\delta_1}{\langle x,x_0\rangle_s}\right)^1 \nonumber +\frac{(\lambda+1)(\lambda+2)}{2!}\left(\frac{-\delta_1}{\langle x,x_0\rangle_s}\right)^2+\cdots\\&+&\frac{(\lambda+1)_{n-1}}{(n-1)!}\left(\frac{-\delta_1}{\langle \nonumber x,x_0\rangle_s}\right)^{n-1}+\frac{(\lambda+1)_{n}}{(n)!}\left(\frac{1}{1+\xi_4}\right)^{\lambda+n+1}\left(\frac{-\delta_1}{\langle x,x_0\rangle_s}\right)^{n} \Bigg],
\end{eqnarray}

\begin{eqnarray}\label{34}
E&=&\frac{x_0(s-1)}{\langle x,x_0\rangle_s^{\lambda+1} \left(1+\frac{\delta_1}{\langle x,x_0\rangle_s}\right)^{\lambda+1}}
\\&+&\frac{x_0(s-1)}{\langle x,x_0\rangle_s^{\lambda+1}}\Bigg[1+\frac{\lambda+1}{1}\left(\frac{-\delta_1}{\langle x,x_0\rangle_s}\right)^1 \nonumber +\frac{(\lambda+1)(\lambda+2)}{2!}\left(\frac{-\delta_1}{\langle x,x_0\rangle_s}\right)^2+\cdots\\&+&\frac{(\lambda+1)_{n}}{(n)!}\left(\frac{-\delta_1}{\langle \nonumber x,x_0\rangle_s}\right)^{n}+\frac{(\lambda+1)_{n+1}}{(n+1)!}\left(\frac{1}{1+\xi_5}\right)^{\lambda+n+2}\left(\frac{-\delta_1}{\langle x,x_0\rangle_s}\right)^{n+1} \Bigg].\nonumber
\end{eqnarray}

By  (\ref{1}),  $\xi_3,  \xi_4, \xi_5 \in[-3/4,3/4]$.
Thus:
\begin{eqnarray}\label{15}
\left(\frac{1}{1+\xi_i}\right)\leq \left(\frac{1}{1-3/4}\right)\leq 4  \ \hbox{for}\  i=3,4,5.
\end{eqnarray}
Thus  Formulas\,(\ref{1}, \ref{15}, \ref{31}, \ref{32}, \ref{33}, \ref{34}, \ref{8q}, \ref{8qq}, \ref{taylor}), together with the vanishing property of $a(t)$, we  obtain:
for $x\in \left[-2x_0, 0\right]\bigcap \left(I(x_0,4\delta_0)\bigcup I(-x_0,4\delta_0)\right)^c $
\begin{eqnarray}\label{poisson inequality 2}
\left|\SH_\lambda a(x)\right|&\lesssim&\int\int_{-1}^1 |a(t)|\left|\frac{x-x_0s}{\langle x,x_0\rangle_s^{\lambda+1}}\left(\frac{-\delta_1}{\langle x,x_0\rangle_s}\right)^{n+1}\right|(1+s)(1-s^2)^{\lambda-1}ds|t|^{2\lambda}dt\nonumber\\&+&
\int\int_{-1}^1 |a(t)|\left|\frac{x_0-t}{\langle x,x_0\rangle_s^{\lambda+1}}\left(\frac{-\delta_1}{\langle x,x_0\rangle_s}\right)^{n}\right|(1+s)(1-s^2)^{\lambda-1}ds|t|^{2\lambda}dt\nonumber\\
&+&\int\int_{-1}^1 |a(t)|\left|\frac{x_0(s-1)}{\langle x,x_0\rangle_s^{\lambda+1}}\left(\frac{-\delta_1}{\langle x,x_0\rangle_s}\right)^{n+1}\right|(1+s)(1-s^2)^{\lambda-1}ds|t|^{2\lambda}dt\nonumber\\
&\lesssim&  |I(x_{0} ,\delta_{0})|_{\lambda}^{1-(1/p)} \int_{-1}^{1} \frac{(\delta_0)^{n+1}}{\left(\langle x,x_0\rangle_s\right)^{n/2+\lambda+1}} (1+s)(1-s^2)^{\lambda-1}ds \nonumber
\\&+&\int\int_{-1}^1 |a(t)|\left|\frac{x_0(s-1)}{\langle x,x_0\rangle_s^{\lambda+1}}\left(\frac{-\delta_1}{\langle x,x_0\rangle_s}\right)^{n+1}\right|(1+s)(1-s^2)^{\lambda-1}ds|t|^{2\lambda}dt.\nonumber\\
\end{eqnarray}
Notice that $\langle x,x_0\rangle_s=x^2+x_0^2-2xx_0s\geq(1-s^2)x_0^2$ holds for $\forall x\in \RR$. Thus we could obtain the following inequality:
\begin{eqnarray}\label{100}
& &\int\int_{-1}^1 |a(t)|\left|\frac{x_0(s-1)}{\langle x,x_0\rangle_s^{\lambda+1}}\left(\frac{-\delta_1}{\langle x,x_0\rangle_s}\right)^{n+1}\right|(1+s)(1-s^2)^{\lambda-1}ds|t|^{2\lambda}dt
\nonumber\\&\lesssim&
\int\int_{-1}^1 |a(t)|\left|\frac{x_0(s-1)}{x_0\left(1-s^2\right)^{1/2}}\frac{(\delta_0)^{n+1}}{\left(\langle x,x_0\rangle_s\right)^{n/2+\lambda+1}}\right|(1+s)(1-s^2)^{\lambda-1}ds|t|^{2\lambda}dt
\nonumber\\&\lesssim&
|I(x_{0} ,\delta_{0})|_{\lambda}^{1-(1/p)} \int_{-1}^{1} \frac{(\delta_0)^{n+1}}{\left(\langle x,x_0\rangle_s\right)^{n/2+\lambda+1}} (1-s^2)^{\lambda-1/2}ds
\nonumber\\&\lesssim&
|I(x_{0} ,\delta_{0})|_{\lambda}^{1-(1/p)} \int_{-1}^{1} \frac{(\delta_0)^{n+1}}{||x|-|x_0||^n\left(\langle x,x_0\rangle_s\right)^{\lambda+1}} (1-s^2)^{\lambda-1/2}ds.
\end{eqnarray}
By Proposition\,\ref{estimate d}, Formula\,(\ref{100}) implies that:
\begin{eqnarray}\label{102}
& &|I(x_{0} ,\delta_{0})|_{\lambda}^{1-(1/p)} \int_{-1}^{1} \frac{(\delta_0)^{n+1}}{||x|-|x_0||^n\left(\langle x,x_0\rangle_s\right)^{\lambda+1}} (1-s^2)^{\lambda-1/2}ds
\nonumber\\&\lesssim&
|I(x_{0} ,\delta_{0})|_{\lambda}^{1-(1/p)} \frac{(\delta_0)^{n+1}}{\left||x|-|x_0|\right|^{n+2}\left||x|+|x_0|\right|^{2\lambda}}.
\end{eqnarray}
Formulas\,(\ref{101}, \ref{102}) imply the following inequality holds when $x\in \left[-2x_0, 0\right]\bigcap \left(I(x_0,4\delta_0)\bigcup I(-x_0,4\delta_0)\right)^c $:
\begin{eqnarray}\label{103}
\left|\SH_\lambda a(x)\right|\leq C |I(x_{0} ,\delta_{0})|_{\lambda}^{1-(1/p)} \frac{(\delta_0)^{n+1}}{\left||x|-|x_0|\right|^{n+2}\left||x|+|x_0|\right|^{2\lambda}}.
\end{eqnarray}
From Formulas\,(\ref{103}, \ref{104}), we obtain that:
\begin{eqnarray}\label{105}
\left|\SH_\lambda a(x)\right|\leq C |I(x_{0} ,\delta_{0})|_{\lambda}^{1-(1/p)} \frac{(\delta_0)^{n+1}}{\left||x|-|x_0|\right|^{n+2}\left||x|+|x_0|\right|^{2\lambda}}
\end{eqnarray}
holds for $x\in\left(I(x_0,4\delta_0)\bigcup I(-x_0,4\delta_0)\right)^c$.
Because $\frac{2\lambda}{2\lambda+1}< p\leq1$, we could get $ 0\leq 2\lambda(1-p)\leq\frac{2\lambda}{2\lambda+1}<1$. Thus the following holds:
\begin{eqnarray*}
|x|^{2\lambda(1-p)}\leq ||x|-|x_0||^{2\lambda(1-p)}+|x_0|^{2\lambda(1-p)}.
\end{eqnarray*}
Therefore:
\begin{eqnarray}\label{10}
II&\leq&C |I(x_{0} ,\delta_{0})|_{\lambda}^{p-1}\int_{I^c} \Bigg(\frac{(\delta_0)^{n+1}}{\left||x|-|x_0|\right|^{n+2}}\Bigg)^p\bigg(||x|-|x_0||^{2\lambda(1-p)}+|x_0|^{2\lambda(1-p)}\bigg)dx \nonumber
\\&\approx&\nonumber  |I(x_{0} ,\delta_{0})|_{\lambda}^{p-1}\int_{I^c}\Bigg(\frac{(\delta_0)^{n+1}}{\left||x|-|x_0|\right|^{n+2}}\Bigg)^p\bigg(||x|-|x_0||^{2\lambda(1-p)}\bigg)dx
\\&+&\nonumber |I(x_{0} ,\delta_{0})|_{\lambda}^{p-1}\int_{I^c} \Bigg(\frac{(\delta_0)^{n+1}}{\left||x|-|x_0|\right|^{n+2}}\Bigg)^p\bigg(|x_0|^{2\lambda(1-p)}\bigg)dx
\\&=&II_1 + II_2.
\end{eqnarray}
$\kappa$ and p satisfy the relation: $\displaystyle{\kappa=2\left[(2\lambda+1)\frac{1-p}{p}\right] }, n=\kappa/2.$ Therefore we could have:
$$(n+2)p+2\lambda(p-1)>1 \ and \  (n+2)p>1.$$
Then the integral of $II_1$ and $II_2$ converge:
\begin{eqnarray}\label{6}
II_1&=&C |I(x_{0} ,\delta_{0})|_{\lambda}^{p-1}\int_{I^c}\Bigg(\frac{(\delta_0)^{n+1}}{\left||x|-|x_0|\right|^{n+2}}\Bigg)^p\bigg(||x|-|x_0||^{2\lambda(1-p)}\bigg)dx  \nonumber\\ &\leq&
C  |I(x_{0} ,\delta_{0})|_{\lambda}^{p-1}(\delta_0)^{(n+1)p}\int_{I^c} \frac{1}{\left||x|-|x_0|\right|^{(n+2)p+2\lambda(p-1)}}dx
\nonumber\\ &\leq&
C  |I(x_{0} ,\delta_{0})|_{\lambda}^{p-1}(\delta_0)^{(n+1)p}\int_{4\delta_0}^{+\infty} \frac{1}{r^{(n+2)p+2\lambda(p-1)}}dr
\nonumber\\ &\leq& C,
\end{eqnarray}

\begin{eqnarray}\label{7}
II_2&=&C  |I(x_{0} ,\delta_{0})|_{\lambda}^{p-1}\int_{I^c} \Bigg(\frac{(\delta_0)^{n+1}}{\left||x|-|x_0|\right|^{n+2}}\Bigg)^p\bigg(|x_0|^{2\lambda(1-p)}\bigg)dx \nonumber\\ &\leq&
C \  x_0^{2\lambda(p-1)} (\delta_0)^{p-1}|x_0|^{2\lambda(1-p)}(\delta_0)^{(n+1)p}(\delta_0)^{-(n+2)p+1}
\nonumber\\ &\leq&
C.
\end{eqnarray}

Thus from Formulas\,(\ref{8}, \ref{10}, \ref{6}, \ref{7}), the theorem is proved.
\end{proof}

\begin{theorem}\label{theorem 2}
If $a(t)$ is  a $p_{\lambda}$-atom,  with vanishing order $\displaystyle{\kappa\geq2\left[(2\lambda+1)\frac{1-p}{p}\right]}$ then the following holds: $$\frac{2\lambda}{2\lambda+1}< p\leq 1\ \ \ \int_{\RR}|\left(a*_{\lambda}P_y\right)|^{p}(x)|x|^{2\lambda}dx\leq C,$$
C is depend on $\lambda$ and p.
\end{theorem}
\begin{proof}
Assume first that $x_0>0.$ Let $\kappa=2\left[(2\lambda+1)\frac{1-p}{p}\right].$
Thus $\kappa$ is an even integer. Let $n=\kappa/2$.
We could write the above integral as:
\begin{eqnarray*}
\int_{\RR}  |\left(a*_{\lambda}P_y\right)|^{p}(x) |x|^{2\lambda}dx
&=&\int_{I(x_0,4\delta_0)\bigcup I(-x_0,4\delta_0)}  |\left(a*_{\lambda}P_y\right)|^{p}(x) |x|^{2\lambda}dx    \\
&+& \int_{(I(x_0,4\delta_0)\bigcup I(-x_0,4\delta_0))^c}  |\left(a*_{\lambda}P_y\right)|^{p}(x) |x|^{2\lambda}dx \\
&=&III+IV.
\end{eqnarray*}

By\,\cite{ZhongKai Li 3}\,(Theorem 3.8)  and Formula\,(\ref{9}),  we could get
\begin{eqnarray*}
III&=&\int_{I(x_0,4\delta_0)\bigcup I(-x_0,4\delta_0)}  |\left(a*_{\lambda}P_y\right)|^{p}(x) |x|^{2\lambda}dx  \\
&\leq& \left(\int_{I(x_0,4\delta_0)\bigcup I(-x_0,4\delta_0)}  |\left(a*_{\lambda}P_y\right)|^{2}(x) |x|^{2\lambda}dx\right)^{p/2}
\left(\int_{I(x_0,4\delta_0)\bigcup I(-x_0,4\delta_0)}  |x|^{2\lambda}dx\right)^{1-p/2} \\
&\leq& C.
\end{eqnarray*}
Next we estimate the following integer:
\begin{eqnarray*}
IV=\int_{(I(x_0,4\delta_0)\bigcup I(-x_0,4\delta_0))^c}  |\left(a*_{\lambda}P_y\right)|^{p}(x) |x|^{2\lambda}dx.
\end{eqnarray*}
By Proposition\,\ref{Poisson-a}, when $x\in \left(I(x_0,4\delta_0)\bigcup I(-x_0,4\delta_0)\right)^c$ we could write $\left(a*_{\lambda}P_y\right)(x)$ as:
$$\left(a*_{\lambda}P_y\right)(x)=c_{\lambda}\int a(t)(\tau_xP_y)(-t)|t|^{2\lambda}dt .$$
Notice that $t\in supp\,a(t)\subseteq I(x_0, \delta_0)$. By the Taylor expansion of formula $\left(1+\frac{\delta_1}{\langle x,x_0\rangle_{y,s}}\right)^{-\lambda-1}$, we could get the following
\begin{eqnarray}\label{exx3}
\frac{y}{\langle x,t\rangle_{y,s}^{\lambda+1}}&=&\frac{y}{\langle x,x_0\rangle_{y,s}^{\lambda+1} \left(1+\frac{\delta_1}{\langle x,x_0\rangle_{y,s}}\right)^{\lambda+1}}\\&=&\nonumber\frac{y}{\langle x,x_0\rangle_{y,s}^{\lambda+1}}\Bigg[1+\frac{\lambda+1}{1}\left(\frac{-\delta_1}{\langle x,x_0\rangle_{y,s}}\right)^1 +\frac{(\lambda+1)(\lambda+2)}{2!}\left(\frac{-\delta_1}{\langle x,x_0\rangle_{y,s}}\right)^2+\cdots\\&+&\nonumber\frac{(\lambda+1)_{n}}{(n)!}\left(\frac{-\delta_1}{\langle x,x_0\rangle_{y,s}}\right)^{n}+\frac{(\lambda+1)_{n+1}}{(n+1)!}\left(\frac{1}{1+\xi}\right)^{\lambda+n+2}\left(\frac{-\delta_1}{\langle x,x_0\rangle_{y,s}
}\right)^{n+1} \Bigg].
\end{eqnarray}

We could see that:
\begin{eqnarray}\label{u3}
\left| \frac{\delta_1}{\langle x,x_0\rangle_{y,s}} \right| \leq \left| \frac{3|t-x_0|}{\left(\langle x,x_0\rangle_{y,s}\right)^{1/2}} \right|\leq \left| \frac{3|t-x_0|}{\left||x|-|x_0|\right|} \right|\leq \left|\frac{3\delta_0}{4\delta_0} \right|=3/4.
\end{eqnarray}
From  (\ref{u3}), we could have: $\xi\in[-3/4,3/4]$. Thus:
\begin{eqnarray}\label{u4}\left(\frac{1}{1+\xi}\right)\leq \left(\frac{1}{1-3/4}\right)\leq 4, \ \ \frac{y}{\langle x,x_0\rangle_{y,s}}\leq \frac{1}{2\langle x,x_0\rangle^{1/2}_{s}}, and \ \ \frac{1}{\langle x,x_0\rangle_{y,s}}\leq\frac{1}{\langle x,x_0\rangle_{s}}\end{eqnarray}
From  (\ref{u4}),  we could have:
\begin{eqnarray*}
\left|\left(\frac{1}{1+\xi}\right)^{\lambda+n+2}\frac{y}{\left(\langle x,x_0\rangle_{y,s}\right)^{\lambda+1}} \frac{{\delta_1}^{n+1}}{\left(\langle x,x_0\rangle_{y,s}\right)^{n+1}}\right|&\leq& C \frac{y}{\left(\langle x,x_0\rangle_{y,s}\right)^{\lambda+1}} \frac{|t-x_0|^{n+1}}{\left(\langle x,x_0\rangle_{y,s}\right)^{n+1}} \left(\langle x,x_0\rangle_{s}\right)^{\frac{n+1}{2}}\\
&\leq&C\frac{|t-x_0|^{n+1}}{\left(\langle x,x_0\rangle_{s}\right)^{n/2+\lambda+1}}.
\end{eqnarray*}

Since we have
\begin{eqnarray}
|x_0+t-2xs|&\leq&|x_0-xs|+|t-xs|\\\nonumber
&\leq&\left(\langle x,x_0\rangle_s\right)^{1/2}+\left(\langle x,x_0\rangle_s\right)^{1/2}+|t-x_0| \\\nonumber
&\leq&3\left(\langle x,x_0\rangle_s\right)^{1/2},\nonumber
\end{eqnarray}
then the following inequality holds:
\begin{eqnarray}\label{taylor11}
|\delta_1|\leq3|t-x_0|\left(\langle x,x_0\rangle_s\right)^{1/2}.
\end{eqnarray}

Thus for $x\in \left(I(x_0,4\delta_0)\bigcup I(-x_0,4\delta_0)\right)^c$,   by Proposition\,\ref{Poisson-a} and Formula(\ref{exx3}),   we could deduce the following
\begin{eqnarray*}\label{poisson inequality}
\left(a*_{\lambda}P_y\right)(x)&=&c_{\lambda}2^{\lambda+1/2}\pi^{-1}\lambda\Gamma(\lambda+1/2)\int a(t)\int_{-1}^1 \frac{y}{\langle x,x_0\rangle_{y,s}^{\lambda+1}}\Bigg[1+\frac{\lambda+1}{1}\left(\frac{-\delta_1}{\langle x,x_0\rangle_{y,s}}\right)^1 \nonumber\\
&+&\frac{(\lambda+1)(\lambda+2)}{2!}\left(\frac{-\delta_1}{\langle x,x_0\rangle_{y,s}}\right)^2+\cdots+\frac{(\lambda+1)_{n}}{(n)!}\left(\frac{-\delta_1}{\langle x,x_0\rangle_{y,s}}\right)^{n}\nonumber\\
&+&\frac{(\lambda+1)_{n+1}}{(n+1)!}\left(\frac{1}{1+\xi}\right)^{\lambda+n+2}\left(\frac{-\delta_1}{\langle x,x_0\rangle_{y,s}
}\right)^{n+1} \Bigg](1+s)(1-s^2)^{\lambda-1}ds|t|^{2\lambda}dt.
\end{eqnarray*}
Then the above formula together with the vanishing property of a(t) we obtain:
\begin{eqnarray}
\left|\left(a*_{\lambda}P_y\right)(x)\right|\leq C |I(x_{0} ,\delta_{0})|_{\lambda}^{1-(1/p)} \int_{-1}^{1} \frac{(\delta_0)^{n+1}}{\left(\langle x,x_0\rangle_s\right)^{n/2+\lambda+1}} (1+s)(1-s^2)^{\lambda-1}ds.
\end{eqnarray}
By  (\ref{101}),
\begin{eqnarray}\label{106}
\left|\left(a*_{\lambda}P_y\right)(x)\right|\leq C |I(x_{0} ,\delta_{0})|_{\lambda}^{1-(1/p)} \frac{(\delta_0)^{n+1}}{\left||x|-|x_0|\right|^{n+2}\left||x|+|x_0|\right|^{2\lambda}}
\end{eqnarray}
holds for $x\in\left(I(x_0,4\delta_0)\bigcup I(-x_0,4\delta_0)\right)^c$. Thus the Theorem\,\ref{theorem 2} is proved in the same way as Theorem\,\ref{p atom}. This proves the theorem.
\end{proof}

\begin{theorem}\label{theorem 3}
If a(t) is  a $p_{\lambda}$-atom, with vanishing order $\displaystyle{\kappa\geq2\left[(2\lambda+1)\frac{1-p}{p}\right]}$  then $$(\frac{2\lambda}{2\lambda+1}< p\leq 1)\ \ \ \int_{\RR}|\left(a*_{\lambda}Q_y\right)|^{p}(x)|x|^{2\lambda}dx\leq C ,$$
C is depend on $\lambda$ and p.
\end{theorem}
\begin{proof}
Assume first that $x_0>0.$ Let $\kappa=2\left[(2\lambda+1)\frac{1-p}{p}\right].$
Thus $\kappa$ is an even integer. Let $n=\kappa/2$.
We could write the above integral as:
\begin{eqnarray*}
\int_{\RR}  |\left(a*_{\lambda}Q_y\right)|^{p}(x) |x|^{2\lambda}dx
&=&\int_{I(x_0,4\delta_0)\bigcup I(-x_0,4\delta_0)}  |\left(a*_{\lambda}Q_y\right)|^{p}(x) |x|^{2\lambda}dx    \\
&+& \int_{(I(x_0,4\delta_0)\bigcup I(-x_0,4\delta_0))^c}  |\left(a*_{\lambda}Q_y\right)|^{p}(x) |x|^{2\lambda}dx \\
&=&V+VI.
\end{eqnarray*}

We could have the estimation:
\begin{eqnarray*}
V&=&\int_{I(x_0,4\delta_0)\bigcup I(-x_0,4\delta_0)}  |\left(a*_{\lambda}Q_y\right)|^{p}(x) |x|^{2\lambda}dx  \\
&\leq& \left(\int_{I(x_0,4\delta_0)\bigcup I(-x_0,4\delta_0)}  |\left(a*_{\lambda}Q_y\right)|^{2}(x) |x|^{2\lambda}dx\right)^{p/2}
\left(\int_{I(x_0,4\delta_0)\bigcup I(-x_0,4\delta_0)}  |x|^{2\lambda}dx\right)^{1-p/2} \\
&\leq& C \left(\int_{\RR}|a(x)|^2 |x|^{2\lambda}dx\right)^{p/2} \left(|I(x_0,4\delta_0)|_{\lambda}\right)^{1-(p/2)}   2^{1-(p/2)}\\
&\leq& C.
\end{eqnarray*}
Next we estimate the following integer:
\begin{eqnarray*}
VI=\int_{(I(x_0,4\delta_0)\bigcup I(-x_0,4\delta_0))^c}  |\left(a*_{\lambda}Q_y\right)|^{p}(x) |x|^{2\lambda}dx.
\end{eqnarray*}
By Proposition\,\ref{Poisson-a}, when $x\in \left(I(x_0,4\delta_0)\bigcup I(-x_0,4\delta_0)\right)^c$ we could write $\left(a*_{\lambda}Q_y\right)(x)$ as:
$$\left(a*_{\lambda}Q_y\right)(x)=c_{\lambda}\int a(t)(\tau_xQ_y)(-t)|t|^{2\lambda}dt,$$
where $ supp\,a(t)\subseteq I(x_0, \delta_0)$.
Notice that the following holds
\begin{eqnarray}\label{exxx1}
\frac{x-t}{\langle x,t\rangle_{y,s}^{\lambda+1}}&=&\frac{x-x_0}{\langle x,t\rangle_{y,s}^{\lambda+1}}+\frac{x_0-t}{\langle x,t\rangle_{y,s}^{\lambda+1}}
\\&=&\nonumber F+G.
\end{eqnarray}

By the  Taylor expansion of  $\left(1+\frac{\delta_1}{\langle x,x_0\rangle_{y,s}}\right)^{-\lambda-1}, $  when $x\in \left[-2x_0, 0\right]^c\bigcap \left(I(x_0,4\delta_0)\bigcup I(-x_0,4\delta_0)\right)^c $, we could obtain:($ \exists \xi $)

\begin{eqnarray}\label{exxx2}
F&=&\frac{x-x_0}{\langle x,x_0\rangle_{y,s}^{\lambda+1} \left(1+\frac{\delta_1}{\langle x,x_0\rangle_{y,s}}\right)^{\lambda+1}}
\\&=&\nonumber\frac{x-x_0}{\langle x,x_0\rangle_{y,s}^{\lambda+1}}\Bigg[1+\frac{\lambda+1}{1}\left(\frac{-\delta_1}{\langle x,x_0\rangle_{y,s}}\right)^1 \nonumber +\frac{(\lambda+1)(\lambda+2)}{2!}\left(\frac{-\delta}{\langle x,x_0\rangle_{y,s}}\right)^2+\cdots\\&+&\frac{(\lambda+1)_{n}}{(n)!}\left(\frac{-\delta_1}{\langle \nonumber x,x_0\rangle_{y,s}}\right)^{n}+\frac{(\lambda+1)_{n+1}}{(n+1)!}\left(\frac{1}{1+\xi_1}\right)^{\lambda+n+2}\left(\frac{-\delta_1}{\langle x,x_0\rangle_{y,s}}\right)^{n+1} \Bigg]\nonumber,
\end{eqnarray}

and

\begin{eqnarray}\label{exxx3}
G&=&\frac{x_0-t}{\langle x,x_0\rangle_{y,s}^{\lambda+1} \left(1+\frac{\delta_1}{\langle x,x_0\rangle_{y,s}}\right)^{\lambda+1}}
\\&+&\nonumber\frac{x_0-t}{\langle x,x_0\rangle_{y,s}^{\lambda+1}}\Bigg[1+\frac{\lambda+1}{1}\left(\frac{-\delta_1}{\langle x,x_0\rangle_{y,s}}\right)^1 \nonumber +\frac{(\lambda+1)(\lambda+2)}{2!}\left(\frac{-\delta}{\langle x,x_0\rangle_{y,s}}\right)^2+\cdots\\&+&\frac{(\lambda+1)_{n-1}}{(n-1)!}\left(\frac{-\delta_1}{\langle \nonumber x,x_0\rangle_{y,s}}\right)^{n-1}+\frac{(\lambda+1)_{n}}{(n)!}\left(\frac{1}{1+\xi_2}\right)^{\lambda+n+1}\left(\frac{-\delta_1}{\langle x,x_0\rangle_{y,s}}\right)^{n} \Bigg].
\end{eqnarray}

We have to estimate the size of formula $\frac{\delta_1}{\langle x,a\rangle_{y,s}}$.
Since $|\delta_1|\leq3|t-x_0|\left(\langle x,x_0\rangle_{s}\right)^{1/2}$, and $|x-x_0|\thickapprox\left(\langle x,x_0\rangle_s\right)^{1/2}$, when $x\in \left[-2x_0, 0\right]^c\bigcap \left(I(x_0,4\delta_0)\bigcup I(-x_0,4\delta_0)\right)^c $. Then:
$$\left| \frac{\delta_1}{\langle x,x_0\rangle_{y,s}} \right| \leq \left| \frac{3|t-x_0|}{\left(\langle x,x_0\rangle_{y,s}\right)^{1/2}} \right|\leq \left| \frac{3|t-x_0|}{\left||x|-|x_0|\right|} \right|\leq \left|\frac{3\delta_0}{4\delta_0} \right|=3/4,$$
\begin{eqnarray*}
\left(\frac{1}{1+\xi_1}\right)\leq \left(\frac{1}{1-3/4}\right)\leq 4 \ and\  \left(\frac{1}{1+\xi_2}\right)\leq \left(\frac{1}{1-3/4}\right)\leq 4.
\end{eqnarray*}

Let $d\nu(s)$ denote $(1+s)(1-s^2)^{\lambda-1}ds$. Then by Proposition\,\ref{Poisson-a} and Formulas\,(\ref{exxx1}, \ref{exxx2}, \ref{exxx3})   for  $x\in \left[-2x_0, 0\right]^c\bigcap \left(I(x_0,4\delta_0)\bigcup I(-x_0,4\delta_0)\right)^c $, we could get
\begin{eqnarray*}\label{Q poisson inequality}
\left(a*_{\lambda}Q_y\right)(x)&\thickapprox&\int\int_{-1}^1 a(t) \frac{x-x_0}{\langle x,x_0\rangle_{y,s}^{\lambda+1}}\Bigg[1+\frac{\lambda+1}{1}\left(\frac{-\delta_1}{\langle x,x_0\rangle_{y,s}}\right)^1+\frac{(\lambda+1)_2}{2!}\left(\frac{-\delta_1}{\langle x,x_0\rangle_{y,s}}\right)^2+\cdots\nonumber
\\&+&\frac{(\lambda+1)_{n}}{(n)!}\left(\frac{-\delta_1}{\langle x,x_0\rangle_{y,s}}\right)^{n}+\frac{(\lambda+1)_{n+1}}{(n+1)!}\left(\frac{1}{1+\xi_1}\right)^{\lambda+n+2}\left(\frac{-\delta_1}{\langle x,x_0\rangle_{y,s}}\right)^{n+1} \Bigg]d\nu(s)|t|^{2\lambda}dt\nonumber\\&+&\int\int_{-1}^1 a(t)\frac{x_0-t}{\langle x,x_0\rangle_{y,s}^{\lambda+1}}\Bigg[1+\frac{\lambda+1}{1}\left(\frac{-\delta_1}{\langle x,x_0\rangle_{y,s}}\right)^1 \nonumber +\frac{(\lambda+1)_2}{2!}\left(\frac{-\delta_1}{\langle x,x_0\rangle_{y,s}}\right)^2 +\cdots\nonumber\\&+&\frac{(\lambda+1)_{n-1}}{(n-1)!}\left(\frac{-\delta_1}{\langle x,x_0\rangle_{y,s}}\right)^{n-1}+\frac{(\lambda+1)_{n}}{(n)!}\left(\frac{1}{1+\xi_2}\right)^{\lambda+n+1}\left(\frac{-\delta_1}{\langle x,x_0\rangle_{y,s}}\right)^{n} \Bigg]d\nu(s)|t|^{2\lambda}dt.                              \end{eqnarray*}
Notice that $ 0<\langle x,x_0\rangle_s\leq\langle x,x_0\rangle_{y,s} $. Then the above formula together with the vanishing property of a(t),   for  $x\in \left[-2x_0, 0\right]^c\bigcap \left(I(x_0,4\delta_0)\bigcup I(-x_0,4\delta_0)\right)^c $, we  obtain  :
\begin{eqnarray}\label{109}
\left|\left(a*_{\lambda}Q_y\right)(x)\right|
&\lesssim&\int\int_{-1}^1 |a(t)|\left|\frac{x-x_0}{\langle x,x_0\rangle_s^{\lambda+1}}\left(\frac{-\delta_1}{\langle x,x_0\rangle_s}\right)^{n+1}\right|d\nu(s)|t|^{2\lambda}dt\nonumber\\&+&
\int\int_{-1}^1 |a(t)|\left|\frac{x_0-t}{\langle x,x_0\rangle_s^{\lambda+1}}\left(\frac{-\delta_1}{\langle x,x_0\rangle_s}\right)^{n}\right|d\nu(s)|t|^{2\lambda}dt\nonumber\\
&\leq& C |I(x_{0} ,\delta_{0})|_{\lambda}^{1-(1/p)} \int_{-1}^{1} \frac{(\delta_0)^{n+1}}{\left(\langle x,x_0\rangle_s\right)^{n/2+\lambda+1}} (1+s)(1-s^2)^{\lambda-1}ds.
\end{eqnarray}
We could see that:
\begin{eqnarray}\label{bu1}
\frac{x-t}{\langle x,t\rangle_{y,s}^{\lambda+1}}&=&\frac{x-x_0s}{\langle x,t\rangle_{y,s}^{\lambda+1}}+\frac{x_0-t}{\langle x,t\rangle_{y,s}^{\lambda+1}}+\frac{x_0(s-1)}{\langle x,t\rangle_{y,s}^{\lambda+1}}
\\&=&\nonumber C_1+D_1+E_1.
\end{eqnarray}
Then  we need to estimate $\left(a*_{\lambda}Q_y\right)(x)$ when $x\in \left[-2x_0, 0\right]\bigcap \left(I(x_0,4\delta_0)\bigcup I(-x_0,4\delta_0)\right)^c $.
By the Taylor expansion, we could obtain:

\begin{eqnarray}\label{bu2}
C_1&=&\frac{x-x_0s}{\langle x,x_0\rangle_{y,s}^{\lambda+1} \left(1+\frac{\delta_1}{\langle x,x_0\rangle_{y,s}}\right)^{\lambda+1}}
\\&=&\frac{x-x_0s}{\langle x,x_0\rangle_{y,s}^{\lambda+1}}\Bigg[1+\frac{\lambda+1}{1}\left(\frac{-\delta_1}{\langle x,x_0\rangle_{y,s}}\right)^1 \nonumber +\frac{(\lambda+1)(\lambda+2)}{2!}\left(\frac{-\delta_1}{\langle x,x_0\rangle_{y,s}}\right)^2+\cdots\\&+&\frac{(\lambda+1)_{n}}{(n)!}\left(\frac{-\delta_1}{\langle \nonumber x,x_0\rangle_{y,s}}\right)^{n}+\frac{(\lambda+1)_{n+1}}{(n+1)!}\left(\frac{1}{1+\xi_3}\right)^{\lambda+n+2}\left(\frac{-\delta_1}{\langle x,x_0\rangle_{y,s}}\right)^{n+1} \Bigg]\nonumber,
\end{eqnarray}

\begin{eqnarray}\label{bu3}
D_1&=&\frac{x_0-t}{\langle x,x_0\rangle_{y,s}^{\lambda+1} \left(1+\frac{\delta_1}{\langle x,x_0\rangle_{y,s}}\right)^{\lambda+1}}
\\&+&\frac{x_0-t}{\langle x,x_0\rangle_{y,s}^{\lambda+1}}\Bigg[1+\frac{\lambda+1}{1}\left(\frac{-\delta_1}{\langle x,x_0\rangle_{y,s}}\right)^1 \nonumber +\frac{(\lambda+1)(\lambda+2)}{2!}\left(\frac{-\delta_1}{\langle x,x_0\rangle_{y,s}}\right)^2+\cdots\\&+&\frac{(\lambda+1)_{n-1}}{(n-1)!}\left(\frac{-\delta_1}{\langle \nonumber x,x_0\rangle_{y,s}}\right)^{n-1}+\frac{(\lambda+1)_{n}}{(n)!}\left(\frac{1}{1+\xi_4}\right)^{\lambda+n+1}\left(\frac{-\delta_1}{\langle x,x_0\rangle_{y,s}}\right)^{n} \Bigg],
\end{eqnarray}
and
\begin{eqnarray}\label{bu4}
E_1&=&\frac{x_0(s-1)}{\langle x,x_0\rangle_{y,s}^{\lambda+1} \left(1+\frac{\delta_1}{\langle x,x_0\rangle_{y,s}}\right)^{\lambda+1}}
\\&+&\frac{x_0(s-1)}{\langle x,x_0\rangle_{y,s}^{\lambda+1}}\Bigg[1+\frac{\lambda+1}{1}\left(\frac{-\delta_1}{\langle x,x_0\rangle_{y,s}}\right)^1 \nonumber +\frac{(\lambda+1)(\lambda+2)}{2!}\left(\frac{-\delta_1}{\langle x,x_0\rangle_{y,s}}\right)^2+\cdots\\&+&\frac{(\lambda+1)_{n}}{(n)!}\left(\frac{-\delta_1}{\langle \nonumber x,x_0\rangle_{y,s}}\right)^{n}+\frac{(\lambda+1)_{n+1}}{(n+1)!}\left(\frac{1}{1+\xi_5}\right)^{\lambda+n+2}\left(\frac{-\delta_1}{\langle x,x_0\rangle_{y,s}}\right)^{n+1} \Bigg].\nonumber
\end{eqnarray}

Notice that $ 0<\langle x,x_0\rangle_s\leq\langle x,x_0\rangle_{y,s} $. Then by (\ref{1}), we could have: $\xi_3,  \xi_4, \xi_5 \in[-3/4,3/4]$.
Thus:
\begin{eqnarray}\label{108}
\left(\frac{1}{1+\xi_i}\right)\leq \left(\frac{1}{1-3/4}\right)\leq 4  \ \hbox{for}\  i=3,4,5.
\end{eqnarray}
Notice that $ 0<\langle x,x_0\rangle_s\leq\langle x,x_0\rangle_{y,s} $. Thus Formulas\,(\ref{bu1}, \ref{bu2}, \ref{bu3}, \ref{bu4}, \ref{1}, \ref{108}, \ref{8q}, \ref{8qq}, \ref{taylor}), together with the vanishing property of a(t), we  obtain the following inequality:
\begin{eqnarray}\label{poisson inequality 000}
\left|\left(a*_{\lambda}Q_y\right)(x)\right|&\lesssim&\int\int_{-1}^1 |a(t)|\left|\frac{x-x_0s}{\langle x,x_0\rangle_s^{\lambda+1}}\left(\frac{-\delta_1}{\langle x,x_0\rangle_s}\right)^{n+1}\right|(1+s)(1-s^2)^{\lambda-1}ds|t|^{2\lambda}dt\nonumber\\&+&
\int\int_{-1}^1 |a(t)|\left|\frac{x_0-t}{\langle x,x_0\rangle_s^{\lambda+1}}\left(\frac{-\delta_1}{\langle x,x_0\rangle_s}\right)^{n}\right|(1+s)(1-s^2)^{\lambda-1}ds|t|^{2\lambda}dt\nonumber\\
&+&\int\int_{-1}^1 |a(t)|\left|\frac{x_0(s-1)}{\langle x,x_0\rangle_s^{\lambda+1}}\left(\frac{-\delta_1}{\langle x,x_0\rangle_s}\right)^{n+1}\right|(1+s)(1-s^2)^{\lambda-1}ds|t|^{2\lambda}dt\nonumber\\
&\lesssim&  |I(x_{0} ,\delta_{0})|_{\lambda}^{1-(1/p)} \int_{-1}^{1} \frac{(\delta_0)^{n+1}}{\left(\langle x,x_0\rangle_s\right)^{n/2+\lambda+1}} (1+s)(1-s^2)^{\lambda-1}ds \nonumber
\\&+&\int\int_{-1}^1 |a(t)|\left|\frac{x_0(s-1)}{\langle x,x_0\rangle_s^{\lambda+1}}\left(\frac{-\delta_1}{\langle x,x_0\rangle_s}\right)^{n+1}\right|(1+s)(1-s^2)^{\lambda-1}ds|t|^{2\lambda}dt.
\end{eqnarray}
for $x\in \left[-2x_0, 0\right]\bigcap \left(I(x_0,4\delta_0)\bigcup I(-x_0,4\delta_0)\right)^c $.
Similar to Formulas\,(\ref{poisson inequality 2}, \ref{poisson inequality 000}) could imply that: for $x\in \left[-2x_0, 0\right]\bigcap \left(I(x_0,4\delta_0)\bigcup I(-x_0,4\delta_0)\right)^c $
\begin{eqnarray}\label{poisson inequality 001}
\left|\left(a*_{\lambda}Q_y\right)(x)\right|\lesssim|I(x_{0} ,\delta_{0})|_{\lambda}^{1-(1/p)} \frac{(\delta_0)^{n+1}}{\left||x|-|x_0|\right|^{n+2}\left||x|+|x_0|\right|^{2\lambda}}.
\end{eqnarray}
Formulas\,(\ref{101}, \ref{109}) imply that: for $x\in \left[-2x_0, 0\right]^c\bigcap \left(I(x_0,4\delta_0)\bigcup I(-x_0,4\delta_0)\right)^c $
\begin{eqnarray}\label{poisson inequality 002}
\left|\left(a*_{\lambda}Q_y\right)(x)\right|\lesssim|I(x_{0} ,\delta_{0})|_{\lambda}^{1-(1/p)} \frac{(\delta_0)^{n+1}}{\left||x|-|x_0|\right|^{n+2}\left||x|+|x_0|\right|^{2\lambda}}.
\end{eqnarray}
Thus similar to (\ref{105}), we could obtain that
\begin{eqnarray*}
VI=\int_{(I(x_0,4\delta_0)\bigcup I(-x_0,4\delta_0))^c}  |\left(a*_{\lambda}Q_y\right)|^{p}(x) |x|^{2\lambda}dx\leq C.
\end{eqnarray*}
This proves the theorem.
\end{proof}
In the way similar to Theorem\,\ref{theorem 2} and Theorem\,\ref{theorem 3}, we could obtain the following Proposition
\begin{proposition}\label{theorem 4}
If $a(t)$ is a  $p_{\lambda}$-atom,  then $$\left(\int_{\RR}\left|\partial_y\bigg(P(a_{k})(x, y)\bigg)\right|^{p}|x|^{2\lambda}dx\right)^{1/p}\lesssim \frac{1}{y},\ \ \left(\int_{\RR}\left|\partial_y\bigg(Q(a_{k})(x, y)\bigg)\right|^{p}|x|^{2\lambda}dx\right)^{1/p}\lesssim \frac{1}{y},$$
where $C$ is a  constant depend on $\lambda$ and p.
\end{proposition}

\section{Main result }

\begin{proposition}\label{analy}
Let $\{f_k\}_k$ $(k\in\NN)$ be a sequence of $\lambda$-analytic functions on the set $S\bigcap\{y>0\}$. If $\sum_k |\lambda_k|| f_k|$, $\sum_k |\lambda_k || \partial_x f_k|$  and $\sum_k |\lambda_k|| \partial_y f_k|$ $(k\in\NN)$ converges uniformly on the set $S\bigcap\{y>0\}$, then  $\sum_k \lambda_k f_k$ $(k\in\NN)$ is a $\lambda$-analytic function on the set $S\bigcap\{y>0\}$.
\end{proposition}
\begin{proof}
We denote $D_{\bar{z}}$ as $$D_{\bar{z}}=\frac{1}{2}(D_{x}+i\partial_y).$$
 A function F(z)=F(x,y)=u(x,y)+iv(x,y) is a $\lambda$-analytic function if and only if F(z)  satisfies the following  $\lambda$-Cauchy-Riemann equations:
$$\left\{\begin{array}{ll}
                                    D_xu-\partial_y v=0,&  \\
                                    \partial_y u +D_xv=0.&
                                 \end{array}\right.
$$
Then the $\lambda$-Cauchy-Riemann equations could be replaced by $D_{\bar{z}}F(z)=0$.
Thus a function F(z)=F(x,y)=u(x,y)+iv(x,y) is a $\lambda$-analytic function if and only if $$D_{\bar{z}}F(z)=0.$$

Notice that $\sum_k |\lambda_k|| f_k|$, $\sum_k |\lambda_k ||\partial_x f_k|$ and $\sum_k |\lambda_k ||\partial_y f_k|$ $(k\in\NN)$ converge uniformly on the set $S\bigcap\{y>0\}$. Thus we could have
$$ \sum_k\lambda_k \partial_y f_k=\partial_y\left(\sum_k \lambda_k  f_k\right),\ \ \ \sum_k\lambda_k \partial_x f_k=\partial_x\left(\sum_k \lambda_k  f_k\right)$$
on the set $S\bigcap\{y>0\}$.
Thus $\sum_k \lambda_k f_k$ $(k\in\NN)$ is a $\lambda$-analytic function on the set $S\bigcap\{y>0\}$.
\end{proof}
\begin{lemma}\label{cesa2}
For $0<p\leq1$, $i\in\NN$, if $\sum_i|a_i|^p<\infty$, then
$$\left(\sum_i|a_i|\right)^p\leq \sum_i|a_i|^p .$$
\end{lemma}

\begin{proof}
Without loss of generality, we may assume that
$$\sum_i|a_i|^p=1.$$
In this case, we have $|a_i|\leq1$ for every i, so
$$\left(\sum_i|a_i|\right)\leq \sum_i|a_i|^p|a_i|^{1-p} \leq \sum_i|a_i|^p=1 .$$
Then we could obtain:
$$\left(\sum_i|a_i|\right)^p\leq \sum_i|a_i|^p .$$
This proves the Lemma.
\end{proof}

Let $$ \ f=\sum_{k}\lambda_{k}a_{k}(x)\in H_{\lambda}^{p}(\RR).$$
By Proposition\,\ref{Poisson-a}, $Pa(x, y)+iQa(x, y)$ is a $\lambda$-analytic function. Together with Theorem\,\ref{theorem 3} and Theorem\,\ref{theorem 2}, $P(a_{k})(x, y)+iQ(a_{k})(x, y) \in H_{\lambda}^p(\RR_+^2)$. From Proposition\,\ref{p5}, we could deduce that
$$\sup_{x\in\RR}|P(a_{k})(x, y)+iQ(a_{k})(x, y)| \leq cy^{-(1/p)(1+2\lambda)}.$$
Then by Lemma\,\ref{cesa2}, we could deduce the following for $y>0$
$$\sup_{x\in\RR}\left|\sum_k\lambda_{k}\bigg(P(a_{k})(x, y)+iQ(a_{k})(x, y)\bigg)\right| \leq cy^{-(1/p)(1+2\lambda)}\left(\sum_k|\lambda_k|^p\right)^{1/p}.$$
Notice that the following holds
$$\left|\partial_y\bigg(P(a_{k})(x, y)\bigg)\right|= \left| D_x\bigg(Q(a_{k})(x, y)\bigg)\right|,\ \ \ \left|\partial_y\bigg(Q(a_{k})(x, y)\bigg)\right|= \left| D_x\bigg(P(a_{k})(x, y)\bigg)\right|.$$
Then by  Proposition\,\ref{analy} and Proposition\,\ref{theorem 4}, we could deduce that
$\sum_k\lambda_{k}\bigg(P(a_{k})(x, y)+iQ(a_{k})(x, y)\bigg)\in H_{\lambda}^p(\RR_+^2).$ Thus for $f(x)\in H_{\lambda}^{p}(\RR)$, we could write $Pf(x,y)$, $Qf(x,y)$ as
$$Pf(x,y)=\sum_k\lambda_{k}P(a_{k})(x, y),\ \ \ Qf(x,y)=\sum_k\lambda_{k}Q(a_{k})(x, y). $$
Then by  Theorem\,\ref{theorem 2} and Theorem\,\ref{theorem 3},  we could obtain
\begin{eqnarray*}
\sup_{y>0}\left(\int_{\RR}\left|\sum_{k}\lambda_{k}a_{k}*_{\lambda}P_y (x)+i\sum_{k}\lambda_{k}a_{k}*_{\lambda}Q_y (x)\right|^p|x|^{2\lambda}dx\right)
 \lesssim \sum_{k}|\lambda_{k}|^{p}.
\end{eqnarray*}
Thus we could obtain the following Proposition:
\begin{proposition}\label{hardy 1}
 For any $f\in H_{\lambda}^{p}(\RR)$,   $\frac{2\lambda}{2\lambda+1}< p\leq1,$ we could deduce that $Pf(x, y)+iQf(x, y)\in H_{\lambda}^p(\RR_+^2)$ and the following holds:
 $$ \sup_{y>0}\left(\int_{\RR}|f*_{\lambda}P_y (x)+if*_{\lambda}Q_y (x)|^p|x|^{2\lambda}dx\right) \leq C\|f\|_{H_{\lambda}^{p}(\RR)} .$$
 \end{proposition}

By Proposition\,\ref{hardy 1} and Theorem\,\ref{us1}, we could obtain the main result in this paper
\begin{theorem}\label{tar}
 For any $f\in H_{\lambda}^{p}(\RR)$,   $\frac{2\lambda}{2\lambda+1}< p\leq1,$ we could deduce  the following
 \begin{eqnarray*}
\int_0^{\infty}|(\SF_{\lambda}f)(\xi)|^k|\xi|^{(2\lambda+1)(k-1-k/p)+2\lambda}d\xi\le
c\|f\|_{H_{\lambda}^{p}(\RR)}^p.
\end{eqnarray*}
When $k=p$, we could obtain
$$\int_0^{\infty}|(\SF_{\lambda}f)(\xi)|^p|\xi|^{(2\lambda+1)(p-2)+2\lambda}d\xi\le
c\|f\|_{H_{\lambda}^{p}(\RR)}^p.$$
 \end{theorem}

\end{document}